\documentclass[aap,MSNbibl,citesort,dvips]{arximspdf}
\usepackage{mathbh}

\doi{10.1214/11-AAP795}
\volume{22}
\issue{3}
\pubyear{2012}
\firstpage{1243}
\lastpage{1265}

\makeatletter

\newtheorem{theo}{Theorem}[section]
\newtheorem{lemma}[theo]{Lemma}
\newtheorem{proposition}[theo]{Proposition}

\newproclaim{definition}[theo]{Definition}
\newproclaim{remark}[theo]{Remark}

\newcommand{\ant}{a_n(\psi,\alpha)}
\newcommand{\Ip}{\mathcal{I}}

\makeatother

\begin{document}
\begin{frontmatter}

\title{Optimal stopping problems for some Markov~processes}
\runtitle{Optimal stopping problems for some Markov processes}

\begin{aug}
\author[A]{\fnms{Mamadou} \snm{Ciss\'e}\ead[label=e1]{cissemk@yahoo.fr}},
\author[B]{\fnms{Pierre} \snm{Patie}\thanksref{t2}\ead[label=e2]{ppatie@ulb.ac.be}}
\and
\author[C]{\fnms{Etienne} \snm{Tanr\'e}\corref{}\ead[label=e3]{Etienne.Tanre@inria.fr}}
\runauthor{M. Ciss\'e, P. Patie and E. Tanr\'e}
\affiliation{ENSAE-S\'en\'egal, Universit\'e Libre de Bruxelles and
INRIA}
\address[A]{M. Ciss\'e\\
ENSAE-S\'en\'egal\\
Libert\'e VI extension, VDN, N 26\\
BP 45512-Dakar Fann, S\'en\'egal\\
France\\
\printead{e1}}
\address[B]{P. Patie\\
D\'epartement de Math\'ematiques\\
Universit\'e Libre de Bruxelles\\
Campus de la Plaine C.P. 210\\
B-1050 Bruxelles\\
Belgique\\
\printead{e2}}
\address[C]{E. Tanr\'e\\
EPI Tosca\\
INRIA Sophia-Antipolis M\'editerran\'ee\\
2004, route des Lucioles BP 93\\
06902 Sophia Antipolis Cedex\\
France\\
\printead{e3}}
\end{aug}

\thankstext{t2}{Supported by the
National Centre of Competence in Research ``Conceptual Issues in
Financial Risk Management'' (NCCR-FINRISK).}

\received{\smonth{4} \syear{2010}}
\revised{\smonth{6} \syear{2011}}

%
\begin{abstract}
In this paper, we solve explicitly the optimal stopping problem with
random discounting and an additive functional as cost of observations
for a regular linear diffusion. We also extend the results to the class
of one-sided regular Feller processes. This generalizes the result of
Beibel and Lerche [\textit{Statist. Sinica} \textbf{7} (1997) 93--108]
and [\textit{Teor. Veroyatn. Primen.} \textbf{45} (2000)
657--669] and Irles and Paulsen [\textit{Sequential Anal.} \textbf{23} (2004)
297--316]. Our approach relies on a combination of techniques borrowed
from potential theory and stochastic calculus. We illustrate our
results by detailing some new examples ranging from linear diffusions
to Markov processes of the spectrally negative type.
\end{abstract}

%

\end{frontmatter}

\section{Introduction}
Consider a one-dimensional regular diffusion $X =\break (X_t )_{ t\geq0}$
with state space
$E = (l,r)$, an interval of
$\mathbb{R}$, defined on a filtered probability space
$(\Omega,(\mathcal{F}_t)_{t\geq0},\mathbb{P})$.
We denote by $(\mathbb{P}_x)_{x\in E}$ the family of probability
measures associated to
the process $X$ such that $\mathbb{P}_x(X_0=x)=1$, and by $\mathbb{E}_x$
the associated expectation operator.
Next, let
$\Sigma_{\infty}^X$ be the family of all stopping times with respect to
the filtration
$\mathbb{F} (=(\mathcal{F}_t)_{t \geq0})$.
In this paper, we are first concerned with the study
of the following optimal stopping problem, for any $x\in E$,
%
\begin{equation}\label{eqoptstop}
\sup_{T \in\Sigma_{\infty}^X}
\mathbb{E}_x[e^{-A_T}g(X_T)-C_T],
\end{equation}
where $g$ is a nonnegative continuous function on $E$, $A =
(A_t)_{t\geq0}$ is a continuous additive functional of the form
%
\begin{equation}\label{addfunctional}
A_t=\int_0^t a(X_s)\,ds
\end{equation}
with $a$ a continuous function on $E$ such that $a(x) > 0$ for all
$x\in E$ and, for any $t\geq0$,
%
\begin{equation}\label{eqcostobserv}
C_t = \int_0^t c(X_s)e^{-A_s}\,ds
\end{equation}
with $c$ a nonnegative continuous function on $E$.
We also aim to solve this problem in the case where $X$ is a Feller
process of the spectrally negative type,
that is, when it has only negative jumps. It is of common use to name~$g$
as the reward function,
$C=(C_t)_{t\geq0}$ as the cost of observations
and $A$ as the random discount factor. We mention that in the case $A
\equiv0$ and $c \equiv0$
the problem (\ref{eqoptstop}) has been studied by Dynkin~\cite{Dynk}
and Shiryaev
\cite{shiryayevbook}
in the general framework of Markov processes.
Moreover, the case $c\equiv0$, that is, without cost of observations,
has been intensively studied
in the literature for one-dimensional diffusions. In particular, Salminen
\cite{Salminen} by means
of the Martin boundary, suggested
a solution to this problem in terms of the excessive majorant function.
More recently, Beibel and Lerche
\cite{Beibel-Lerche-97,Beibel-Lerche-00}, relying on martingales
arguments, solve this optimal stopping problem explicitly.
We also mention that, by using standard fluctuation theory, Kyprianou
and Pistorius
\cite{kyprianoupistorius}
offer solution to some optimal stopping problem arising in financial
mathematics for appropriate
diffusions.
A related result on optimal stopping problems for one-dimensional
diffusions with discounting
has been presented by Dayanik and Karatzas~\cite{DayKar}.
They characterize excessive functions via generalized
concavity and determined the value function as the smallest
concave majorant of the reward
function.
In the former case, the value function is given as the solution of a
free boundary value problem associated
to a second order differential operator which is the infinitesimal
generator of the
one-dimensional diffusion $X$. The term \textit{free}, which comes from
the a priori,
unknown region where the problem is investigated, forces one to set up
an artificial
boundary condition of Neumann type to get a well-posed problem. This is
the so-called
smooth fit principle. All these techniques are well explained in the
book of Peskir and Shiryaev
\cite{PeskirShiryaev}.
The literature regarding optimal stopping problems associated to
diffusion with jumps
is more sparse and focused essentially on the study of specific
examples. In this vein,
we mention the paper of Alili and Kyprianou
\cite{Alili-Kyprianou-05} where the authors deal
with the issue of pricing perpetual American put options in a~market
driven by L\'evy processes.
We also indicate that Baurdoux~\cite{Baurdoux-07} solved an optimal
stopping problem
associated to generalized Ornstein--Uhlenbeck processes of the
spectrally negative type.

In this paper, we propose to solve the optimal stopping problem (\ref
{eqoptstop})
with a cost of observations of the form (\ref{eqcostobserv})
for one-dimensional regular diffusions. Our strategy can be described
as follows.
First, by a time change device, we reduce the optimal stopping problem
with random
discounting to a one with a deterministic discount factor but
associated to an appropriate
time change diffusion. Then, by an argument of potential theory, we
transform the problem (\ref{eqoptstop})
to an optimal stopping problem of the same form with a new reward
function but without cost of observations.
We proceed by using a result of Shiryaev~\cite{shiryayevbook} which
states that in our context the optimal
stopping time is the first exit time of the process from a Borel set.
Finally, with this information at hand,
we can use a Doob's $h$-transform technique, with a proper choice of
the excessive function, to transform
our problem to an optimization problem which has been studied in detail
by Beibel and Lerche~\cite{Beibel-Lerche-00}.

The remaining part of the paper is organized as follows. In the next
section, we overview
the basic facts about one-dimensional diffusions. In Section~\ref
{sec3}, we state and prove our main result which
consists of solving the general optimal stopping problem (\ref
{eqoptstop}). We also show, in Section~\ref{sec4},
how to generalize our result to the class of regular one-sided Feller
processes. The last section is
devoted to the treatment of new examples
ranging from linear diffusions to processes with one sided jumps,
such as spectrally negative L\'evy processes and self-similar positive
Markov processes of the spectrally negative type.

\section{Preliminaries}

In this part, we provide some well-known facts about linear diffusions
which can be found, for instance, in It\^o and McKean
\cite{ItoMckean} and in Borodin and
Salminen~\cite{Borodin-Salminen}. We recall that $(\Omega, (\mathcal
{F}_t)_{t\geq0}, \mathbb{P})$
is a filtered probability space. We consider a linear diffusion
$X=(X_t)_{t\geq0}$ with state space
$E$, as the solution to the stochastic differential equation (SDE)
%
\begin{equation} \label{eqsde}
dX_t = b(X_t) \,dt + \sigma(X_t)\,dW_t,
\end{equation}
where $W=(W_t)_{t\geq0 }$ is a one-dimensional Brownian motion.
It is supposed that $\sigma$ and $b$ are continuous and $\sigma(x)>0$
for all $x\in E$.
We assume that $X$ is regular, that is,
\[
\mathbb{P}_x( T_y < +\infty) >0\qquad \forall x,y \in E,
\]
where $T_y= \inf\{ t > 0; X_t=y\}$.
The transition semigroup $(P_t)_{t\geq0}$ maps $C_{b}(E)$, the space of
bounded and continuous functions on $E$, into itself. It follows that $X$
is a Feller process. Moreover, for every $t\geq0$ and $ x\in E$ the
corresponding
measure $A \longmapsto P_t(x,A)$, with
$A$ a Borel set, is absolutely continuous with respect to the speed
measure $m$,
a positive $\sigma$-finite measure on $E$. More specifically, we have
\[
P_t(x,A)= \int_A p_t(x,y) m(dy),
\]
where $p_t(\cdot,\cdot)$ stands for the transition probability density which
may be taken to be positive, jointly
continuous in all variables and symmetric. The scale function $s$ of
$X$ is an increasing continuous function from $E$ to $\mathbb{R}$,
satisfying, for any $a\leq x \leq b$,
\[
\mathbb{P}_x (T_a < T_b)= \frac{s(b) -s(x)}{s(b) -s(a)}
\]
and is given by
\[
s^\prime(x) = \exp\biggl\{- 2 \int^x \frac{b(z)}{\sigma^2(z)} \,dz \biggr\}.
\]
We also recall that the infinitesimal generator $\mathbf{L}$ of $X$ is
the second order differential operator given,
for a function $f\in C^\infty_c(E)$, the space of infinitely
continuously differentiable functions with compact support, by
\[
\mathbf{L} f(x)= \tfrac{1}{2} \sigma^2 (x) f''(x) + b(x) f'(x).
\]
Next, from the general theory of one-dimensional diffusion (see
\cite{ItoMckean}, page~128), the Laplace transform of
the first hitting time $T_y$ is expressed, for any $q>0$, as
%
\begin{equation}\label{defnitfi}
\mathbb{E}_x [ e ^{-q T_y} ] = \cases{
\dfrac{h^+_q(x)}{h^+_q(y)}, &\quad $x\leq y$,\vspace*{2pt}\cr
\dfrac{h^-_q(x)}{h^-_q(y)}, &\quad $x\geq y$,}
\end{equation}
where $h_q^+$ (resp., $h_q^-$) is the increasing (resp., decreasing)
continuous solution to the
differential equation
%
\begin{equation}\label{Sturmliouville}
\mathbf{L} h(x)= q h(x),\qquad x \in E,
\end{equation}
with appropriate conditions at the nonsingular boundary points.
These functions, $h^+_q$ and $h^-_q$, are called the
fundamental solutions of the equation (\ref{Sturmliouville}).
They are linearly independent and their Wronskian is defined by
\[
W_{h^+_q,h^-_q}(x)= h^-_q(x)\,\frac{d }{dx}h^+_q(x) -h^+_q(x) \,\frac
{d}{dx}h^-_q(x),
\]
and the Wronskian with respect to $d / ds(x)$ denoted by $w_q$ is the
constant given by
%
\begin{equation}\label{eqwronskien}
w_q = \frac{W_{h^+_q,h^-_q}(x)}{s'(x)}.
\end{equation}
Moreover, for any $q > 0$, the Green function or the $q$-potential
density $u^q$ is defined as
the Laplace transform of the transition probability density,
that~is,\looseness=-1
\[
u^q(x,y)= \int_0^{+\infty} e ^{-qt} p_t(x,y) \,dt\qquad
\forall x,y \in E .
\]\looseness=0
In particular, we have
%
\begin{equation}\label{qpotentfi}
u^q(x,y)= \cases{
w_q ^{-1} h^+_q(x) h^-_q(y), &\quad $x \leq y$, \cr
w_q ^{-1} h^+_q(y) h^-_q(x), &\quad $x \geq y$.}\vadjust{\goodbreak}
\end{equation}
We say that the process $X$ is recurrent if and only if
$ \lim_ {q \to0} u^q(x,y)= \infty$, for all $x,y \in E$, which is
equivalent to $\mathbb{P}_x(T_y < \infty) =
1$,
for all $x,y\in E$.
A~diffusion which is not recurrent is called transient.
In this case, the potential~$u$ defined as
\[
u(x,y)= \lim_{q \to0} u^q(x,y)
\]
is finite for all $x,y \in E$.
Finally, we mention that if $X$ is transient with $ \lim_{t \to\infty
} X_t = r$ then for any $x\leq y$
\begin{eqnarray*}
u(x,y)& = &\int_0^{+\infty} p_t(x,y) \,dt, \\
& = &s(r)- s(y).
\end{eqnarray*}

\section{Optimal stopping problem for linear diffusions}\label{sec3}
Our aim is now to find the value of the function $\mathcal{V}^A_{g,c}$
defined as the solution of the optimal
stopping problem (\ref{eqoptstop}), that is,
\[
\mathcal{V}^A_{g,c}(x) = \sup_{T \in\Sigma_\infty^X}
\mathbb{E}_x[e^{-A_T}g(X_T)-C_T],
\]
where $A=(A_t)_{t\geq0}$ is a continuous additive functional of the
form (\ref{addfunctional}) and
the cost of observations $C_t = \int_0^t
c(X_s)e^{-A_s}\,ds$ with $c$ and $g$ nonnegative continuous functions.
Our main result is stated in Theorem~\ref{thmmain1} below. It consists
of reducing the
optimal stopping problem (\ref{eqoptstop}) into a new one which can be
described as follows.
On the one hand, it has a modified reward function but without both
cost of observations
and discounting factor. On the other hand, it is associated to a new
diffusion obtained from the original one
by a~random time change and by a~Doob's $h$-transform.
It turns out that solving this latter optimal stopping problem amounts
to finding the solution of an optimization
problem which has been studied by Beibel and Lerche~\cite{Beibel-Lerche-00}.
More precisely, our approach can then be split into the following three steps.
\begin{longlist}[(1)]
\item[(1)] First, we time change $X$ by the inverse of the continuous
increasing functional $A$ and we use the well-known fact that in our
context the two processes have identical hitting time distributions.
Hence, we may consider without loss of generality the problem
(\ref{eqoptstop}) with linear discounting, that is, $A_t = qt$ for some
constant $q>0$.
\item[(2)] Then, we characterize the potential associated to the
functional $C$ and we show how to reduce our problem to an optimal
stopping problem without cost of observations but involving a new
reward function.
\item[(3)] Finally, we borrow an idea of Williams~\cite{Williams74} and
Pitman and Yor~\cite{Pitman-Yor-81} for constructing conditioned
diffusions by the method of $h$-transform. We transform the problem
described in item (2) to an optimal stopping problem without
discounting factor which has been solved by Shiryaev
\cite{shiryayevbook}.
\end{longlist}

\subsection{Time change for nonnegative additive functional}
We start our program by considering that $A$ is a nonnegative
continuous additive
functional of $X$ of the form (\ref{addfunctional})
and we assume, without loss of generality, that $A_{\infty}=\infty$
$\mathbb{P}_x$ a.s. Since $A$ is a
continuous increasing function,
it admits an inverse functional which we denote by $V$ and given, for
all $t\geq0$, by
\begin{eqnarray*}
V_t &=& \inf\{s\geq0; A_s >t\} ,\\
&=& \int_0^t\frac{1}{a(Y_s)}\,ds,
\end{eqnarray*}
where $Y_t =X_{V_t}$ for any $t\geq0 $.
Moreover, if $X$ is the solution to the SDE (\ref{eqsde}), then it
follows, from the It\^o's formula, that $Y$ is the unique
solution to the SDE
\[
dY_t = \biggl(\frac{b}{a}\biggr)(Y_t)\,dt +
\biggl(\frac{\sigma}{\sqrt{a}}\biggr)(Y_t)\,d\hat{W}_t,
\]
where $\hat{W}$ is a Brownian motion with respect to the filtration
$\mathbb{F}^Y = (\mathcal{F}_{V_t})_{t\geq0}$.
The process $Y$ remains a linear diffusion with respect to the
filtration $\mathbb{F}^Y$. In particular, it is a Feller process (see,
e.g., Lamperti
\cite{Lamperti-67}) and its
infinitesimal generator $\mathbf{L}^Y$ takes the form
%
\begin{equation}\label{generatorY}
\mathbf{L}^Yf(x) = \frac{1}{a(x)}\mathbf{L} f(x)
\end{equation}
for a smooth function $f$ on $E$.
Next, we consider an open interval $B\subset E$ and denote by $T^Y_B$
the first exit time
of the process $Y$ from $B$, it is plain that we have the following identity:
\[
T^Y_B = A_{T^X_B} \qquad\mbox{a.s.}
\]
We are now ready to state the following.
\begin{lemma}\label{timechangelema}
For any $x\in E$, we have, with the obvious notation,
%
\begin{eqnarray} \label{eqopt-pbm2}
&&\sup_{T \in\Sigma^X_{\infty}}
\mathbb{E}_x[e^{-A_T}g(X_T)-C_T]\nonumber\\[-8pt]\\[-8pt]
&&\qquad = \sup_{T \in\Sigma^Y_{\infty}}
\mathbb{E}_x\biggl[e^{-T}g(Y_{T})-\int_0^{T}\biggl(\frac{c}{a}\biggr)(Y_s)
e^{-s}\,ds\biggr],\nonumber
\end{eqnarray}
where $Y$ is characterized by its infinitesimal generator (\ref{generatorY}).
\end{lemma}
\begin{pf}
From~\cite{Rogers-Williams-1}, Section III.21, page 277, we have that for every $\mathbb
{F}^X$-stopping time $T$, $A_T$
is an $\mathbb{F}^Y$-stopping time. Thus, we obtain that
\[
\sup_{T \in\Sigma_\infty^X}
\mathbb{E}_x[e^{-A_T}g(X_T)-C_T] = \sup_{S \in
\Sigma_\infty^Y} \mathbb{E}_x\biggl[e^{-S}g(X_{V_S})-\int_0^{V_S}c(X_v)
e^{-A_v}\,dv\biggr].
\]
The proof follows by performing the change of variable $u=A_v$ in the
integral on the right-hand side of the previous identity.\vadjust{\goodbreak}
\end{pf}
Consequently, in the sequel we can assume, without loss of generality, that
the additive functional $A$ is
linear, that is, $A_t = qt$ for some $q>0$ and the cost of observations is
$C_t = \int_0^t c(X_s)e^{-qs}\,ds$.

\subsection{Get rid of the cost of observations}
Let us now introduce the function $\delta$ defined, for any $x\in E$, by
%
\begin{equation}\label{potentiel}
\delta(x) = \mathbb{E}_x[C_\infty].
\end{equation}
In the following, we provide an expression of $\delta$ in terms of the
characteristics
of $X$ and give some conditions under which it is continuous and finite.
\begin{lemma}\label{finidelta}
For any $q>0$ and $x\in E$, we have
\[
\delta(x)= w_q^{-1} \biggl(h^-_q(x)\int_l^{x} h^+_q(y) c(y) m(dy) +
h^+_q(x)\int_x^{r} h^-_q(y) c(y) m(dy)\biggr),
\]
where $w_q$ stands for the Wronskian of $h^-_q$ and $h^+_q$ with
respect to the scale function~$s$, as defined in (\ref {eqwronskien})
and $m$ is the speed measure of $X$. Moreover, if~$c$ satisfies the
integrability condition for any $x\in E$,
%
\begin{equation} \label{33}
\int_l^{r}h^-_q(x\vee y) h^+_q(x\wedge y) c(y) m(dy) < \infty,
\end{equation}
then $\delta$ is continuous and finite on $E$.
\end{lemma}
\begin{pf}
By using Fubini's theorem, we obtain that
\begin{eqnarray*}
\delta(x) & = & \int_0^{\infty}e^{-qs}\mathbb{E}_x[c(X_s)]\,ds = \int_E
u^q(x,y) c(y) m(dy)\\
& = &w_q^{-1} \int_E h^-_q(x\vee y) h^+_q(x\wedge y) c(y) m(dy),
\end{eqnarray*}
where we have used the identity (\ref{qpotentfi}).
The proof of the claims follows readily.
\end{pf}
\begin{remark}
$\!\!$We note that if $q=0$ and $X$ is transient with \mbox{$ \lim_{t \to\infty}
X_t = r$},
then $\delta$ is given by
%
\begin{equation}\label{35}
\delta(x) = \int_E \bigl(s(r)-s(y)\bigr) c(y) m(dy).
\end{equation}
In this case, we mention that Khoshnevisan, Salminen and Yor
\cite{Khosh-Salm-Yor-06}
identify the law of the perpetual integral functional $C_{\infty}$
of a transient diffusion as the law of the first
hitting time of a random time change diffusion.
\end{remark}

We are now ready to state the following.
\begin{lemma}\label{ramenbabel}
If $\delta$ is finite on $E$ then, for any $x\in E$, we have
\[
\sup_{T \in\Sigma_\infty^X}
\mathbb{E}_x[e^{-qT}g(X_T)-C_T]= \sup_{T \in\Sigma_\infty^X}
\mathbb{E}_x\bigl[e^{-qT}\bigl(g(X_T)+\delta(X_T)\bigr)\bigr]-\delta(x).\vadjust{\goodbreak}
\]
\end{lemma}
\begin{pf}
Note that for any $\mathbb{F}$-stopping time $T$, we have the identity
in~law
\[
C_\infty\stackrel{(d)}{=} C_T +
e^{-qT}C_\infty\circ\theta_T,
\]
where $(\theta_t)_{t\geq0}$ stands for the shift operator, that is,
for any $t,s \geq0,
\theta_t w(s)=w(t+s)$.
Since $\delta$ is finite, the strong Markov property yields
\[
\sup_{T \in\Sigma_\infty^X}
\mathbb{E}_x[e^{-qT}g(X_T)-C_T]= \sup_{T \in\Sigma_\infty^X}
\mathbb{E}_x\bigl[e^{-qT}\bigl(g(X_T)+\delta(X_T)\bigr)\bigr]-\delta(x).\quad
\]
\upqed\end{pf}
A nice consequence of the previous result is that the general optimal
stopping problem
(\ref{eqoptstop}) is equivalent to an optimal stopping problem without
cost of observations.
Before stating our next result, let us introduce a few further
notation. Let $\mathcal{O}$ denote all the
open subsets of $E$ containing the starting point $x$ of $X$.
Let $\Sigma_\mathcal{O}$ be the class of stopping times of the form
$T_B= \inf\{ t > 0; X_t \notin B \}$
where $B\in\mathcal{O}$.
\begin{lemma}\label{propshiryayev}
Let $f$ be a continuous nonnegative function.
Then, for any $ q > 0 $ and $x\in E$,
%
\begin{equation}\label{Eqshirayevchp4}
\sup_{T \in\Sigma_\infty^X}
\mathbb{E}_x[ e ^{-qT} f(X_T)]= \sup_ {T_B\in\Sigma_{\mathcal{O}}}
\mathbb{E}_x[e ^{-q T_B} f(X_{T_B})].
\end{equation}
\end{lemma}
\begin{pf} Let $\hat X$ be defined by
\[
\hat X_t=\cases{
X_t, &\quad if $t < \mathbf{e}_q$,\cr
\partial, &\quad if $t \geq\mathbf{e}_q$,}
\]
where $\mathbf{e}_q$ is an exponential variable of parameter $q>0$
taken independent of $\mathbb{F}$
and $\partial$ is a cemetery state.
Note that $\hat X$ is always transient and clearly with a function $f$
as above and using the
convention $f(\partial)=0$, we have, for any $x\in E$
\[
\mathbb{E} _x [f(\hat X _t)] =\mathbb{E} _x [ e^{-qt}f ( X_t )].
\]
Therefore, there is a one-to-one correspondence between
the excessive functions for $\hat X$ and the
$q$-excessive ones for $X$ (see Definition~\ref{defqexcessive} below).
Moreover,
the $q$-excessive functions
of the Feller process $X$
are lower semi-continuous (\cite{hawkes}, Theorem 2.1). Then, from Shiryaev
\cite{shiryayevbook}, Corollary 3,\break page~129, we obtain that
\[
\sup_ {T \in\Sigma^{\hat X}_{\infty} }\mathbb{E}_x [f(\hat X_T)] =
\sup_ {T_B\in\Sigma_{\mathcal{O}}} \mathbb{E}_x[ f(\hat X_{T_B})].
\]
Hence,
\[
\sup_{T \in\Sigma_\infty^X}\mathbb{E}_x[ e ^{-qT} f(X_T)]=
\sup_{T_B\in\Sigma_{\mathcal{O}}} \mathbb{E}_x[e ^{-q T_B}
f(X_{T_B})] .
\]
\upqed\end{pf}

\subsection{Doob's $h$-transform}\label{sub33}

Our aim in this part is to show how to transform an optimal stopping problem
with discounting factor to an optimal stopping problem without\vadjust{\goodbreak}
discounting. To this end,
we recall some basic facts on excessive functions and Doob's
$h$-transform and we refer to the book of Borodin and Salminen
\cite{Borodin-Salminen}, Section II.5, pages 32--35.
\begin{definition}\label{defqexcessive}
A nonnegative measurable function $h\dvtx E \mapsto\mathbb{R} \cup\{
\infty
\}$ is called $q$-excessive, $q\geq0$,
for the process $X$ if the following two statements hold true: for any
$x\in E$,
\begin{longlist}
\item $e ^{-q t} \mathbb{E}_x [ h(X_t)] \leq h(x), t>0$,
\item $ \lim_{t \searrow0} e ^{-q t} \mathbb{E}_x [h(X_t)] =
h(x)$.
\end{longlist}
A $q$-excessive function is called $q$-invariant if for any $x\in E$
and $t\geq0$ we~have
\[
e ^{-q t} \mathbb{E}_x [ h(X_t)] = h(x).
\]
\end{definition}

A function $h$ is $q$-excessive (resp., $q$-invariant) if and only if
the process
$ e ^{-q t} h(X_t)$ is a positive supermartingale (resp., martingale).
For every $y\in(l,r)$ the functions $x \mapsto u^q(x,y)$, $x \mapsto h^+_q(x)$
and $ x \mapsto h^-_q(x)$ are $q$-excessive. These functions are
\textit{minimal}
in the sense that any other arbitrary nontrivial $q$-excessive
function $h$
can be expressed as a linear combination of them.
\begin{definition} \label{defht}
Let $q \geq0$ and $h$ be a $q$-excessive function. For any $x\in E$
such that $0 < h(x) < +\infty$,
and $t\geq0$,
we define the new probability measure
$\mathbb{P}^{h}_x$ as
\[
d\mathbb{P} ^{h}_x = e ^{-q t} \frac{h(X_t)}{h(x)} \,d\mathbb{P}_x
\qquad\mbox{on } \mathcal{F}_t .
\]
The process $X$ under the probability measure $\mathbb{P}^{h}_x $ is
called the Doob's $h$-transform (or $q$-excessive transform) of $X$.
It is also a regular diffusion process and thus a Feller process.
\end{definition}

Next, we use an idea of Williams~\cite{Williams74} and
Pitman and Yor~\cite{Pitman-Yor-81} for constructing conditioned
diffusions by the method of
$h$-transform by means of the Laplace transform of first passage times.
In fact,
we slightly generalize their methodology by considering as
$q$-excessive function the Laplace transform
of the first exit time of an open set by $X$.
To this end, let $B \in\mathcal{O} $ and we recall that we denote by
$T_B$ the first exit time from $B$ by $X$, that is,
\[
T_{B} = \inf\{t > 0; X_t \notin B \}.
\]
For any $x\in E$, we write the Laplace transform of the stopping time
$T_B$ as
\[
\phi^B(x) = \mathbb{E}_x[e^{-q T_B}].
\]
Thus,\vspace*{1pt} without loss of generality, the continuity of $X$ allows us to
restrict~$\mathcal{O}$
to open intervals $(a,b)$
for some $l\leq a < b \leq r$. It is well known that the function~$\phi
^B$ is solution to the
following Sturm--Liouville boundary value problem
\[
\cases{\mathbf{L} u(x) = q u(x), &\quad $x\in(a,b)$, \cr
u(a)=u(b)=1.}\vadjust{\goodbreak}
\]
In other words, $\phi^B$ can be written as a linear combination of the
fundamental solutions $h^-_q$
and $h^+_q$. Setting
\[
h^B(y) = \frac{\phi^B(y)}{\phi^B(x)},
\]
it is plain that the mapping $h^B$ is a $q$-excessive function for $X$.
Thus, as in the
Definition~\ref{defht}, we define the probability measure $\mathbb
{P}_x^{h^B}$ as
%
\begin{equation}\label{defnproba}
d\mathbb{P}_x^{h^B} = e^{-qT_B} h^B(X_{T_B}) \,d\mathbb{P}_x
\qquad\mbox{on
} \mathcal{F}_{T_B^+}.
\end{equation}
The diffusion $X$ under the family of probability measures
$\mathbb{P}^{h^B}=(\mathbb{P}^{h^B}_x)_{ x\in E}$ is transient with
%
\begin{equation}\label{3hpitmanyor}
\mathbb{P}_x^{h^B} ( X_{\xi^-} \in\partial B)=1,
\end{equation}
except if $q=0$, where $\xi$ stands for the lifetime of
$X$ under $\mathbb{P}^{h^B}$. Clearly, the probability $\mathbb
{P}_x^{h^B} ( \xi< \infty)$
is either $1$ or $0$ for all $x$.
Moreover, since $X$ is solution to the SDE (\ref{eqsde}), then under
the probability
$\mathbb{P}^{h^B}$, the diffusion $X$ can be characterized as the
solution of the SDE
\[
dX_t= \bigl( b(X_t) + \log'(h^{B}(X_t) )\sigma^2(X_t) \bigr)\,dt +
\sigma(X_t)\,d\tilde{W}_t,
\]
where $\tilde{W}$ is a standard Brownian motion under the new
probability $\mathbb{P}_x^{h^B}$.
\begin{remark}
As explained in~\cite{Pitman-Yor-81}, we take $\mathbb{P}_x^{h^B}$ to
be defined by the requirement that for each
$x\in B$ the process $X$ runs up to the time $T_B$ has the same law
under $\mathbb{P}_x^{h^B}$
as it does under $\mathbb{P}_x$ conditional on $T_B < \mathbf{e}_q$,
where $\mathbf{e}_q$ is an independent
exponentially distributed random variable with parameter $q>0$.
\end{remark}

We are now ready to state and prove the main theorem of this section.
\begin{theo}\label{thmmain1}
If $\delta$ is finite then solving the problem (\ref{eqoptstop})
amounts to solving the following optimal stopping problem:
%
\begin{equation}\label{eqopt-pbm3}
\sup_{T_B \in\Sigma_{\mathcal{O}}}
\mathbb{E}_x^{h^B}\biggl[\frac{g(X_{T_B})+ \delta(X_{T_B})}{h^B(X_{T_B})} \biggr],
\end{equation}
where the probability $\mathbb{P}_x^{h^B}$ is defined in (\ref{defnproba}).
If there exists an open interval~$B^*$ of $E$ such that
%
\begin{equation}\label{eqopt-pbm4}
\frac{g(u^*)+ \delta(u^*)}{h^{B^*}(u^*)} = \sup_{B \in\mathcal{O},
u\in\partial B} \frac{g(u)+ \delta(u)}{h^B(u)}
\qquad\mbox{where }u^*\in\partial B^*, |u^*|<\infty,\hspace*{-35pt}
\end{equation}
then, the value function of (\ref{eqopt-pbm3}) is given by (\ref
{eqopt-pbm4}) with
the optimal stopping time $T_{B^*}$.
\end{theo}
\begin{remark}
(1) We mention that the optimization problem (\ref{eqopt-pbm4}) has
been studied in detail by Beibel and Lerche~\cite{Beibel-Lerche-00}. We
refer to their paper for more precise information concerning its
solution for all possible choices of the reward function~$g$. We also
point out that, in the specific case $a=0$ and $X$ is a standard
Brownian motion, a similar optimization problem was studied by
Graversen and Peskir~\cite{Graversen-Peskir}.

(2) If $B^*$ is a bounded interval, that is, $B^* =
(u^*_1,u^*_2)$, (\ref{eqopt-pbm4}) reads
\[
\frac{g(u^*_1)+ \delta(u^*_1)}{h^{B^*}(u^*_1)} = \frac{g(u^*_2)+
\delta(u^*_2)}{h^{B^*}(u^*_2)} =
\sup_{B \in\mathcal{O}} \sup_{ u \in\partial B}\frac{g(u)+
\delta
(u)}{h^B(u)}.
\]
\end{remark}
\begin{pf*}{Proof of Theorem~\ref{thmmain1}}
First, we deal with the additive functional~$C$. Since $\delta$ is
finite, we have, from Lemma~\ref{ramenbabel},
that solving the problem (\ref{eqoptstop}) is equivalent to solving the
following
optimal stopping problem without cost of observations:
\[
\sup_{T \in\Sigma_\infty^X}
\mathbb{E}_x\bigl[e^{-qT}\bigl(g(X_T)+\delta(X_T)\bigr)\bigr].
\]
Then, from Lemma~\ref{propshiryayev}, we deduce that
\[
\sup_{T \in\Sigma_\infty^X}
\mathbb{E}_x\bigl[e^{-qT}\bigl(g(X_T)+\delta(X_T)\bigr)\bigr]= \sup_{T_B\in\Sigma
_{\mathcal{O}}}
\mathbb{E}_x\bigl[e ^{-q T_B} \bigl(g(X_{T_B})+\delta(X_{T_B})\bigr)\bigr] .
\]
Next, we use the Doob's $h$-transform device. Let $\mathbb{P}_x^{h^B}$
be the probability measure defined in (\ref{defnproba}),
then we have
\[
\sup_{T_B\in\Sigma_{\mathcal{O}}} \mathbb{E}_x\bigl[e ^{-q T_B}
\bigl(g(X_{T_B})+\delta(X_{T_B})\bigr)\bigr]
= \sup_{T_B\in\Sigma_{\mathcal{O}}} \mathbb{E}_x^{h^B}\biggl[ \frac
{g(X_{T_B})+\delta(X_{T_B})}{h^B(X_{T_B})}\biggr],
\]
which completes the proof of the first assertion.
Finally, since $X$ under $\mathbb{P}_x^{h^B}$ is transient,
it is stated above that $T_B < \infty$ a.s.,
the value function of the last optimal stopping problem is
the solution to the following optimization problem:
\[
\sup_{B \in\mathcal{O}} \sup_{ u \in\partial B}
\frac{g(u)+ \delta(u)}{h^B(u)}.
\]
\upqed\end{pf*}

\section{Extension to one-sided regular Feller processes}\label{sec4}
Let us now consider~$X$ to be the c\`adl\`ag modification of a
one-sided regular
Feller process defined on a filtered probability space
$(\Omega, (\mathcal{F}_t)_{t\geq0}, \mathbb{P})$ and taking values in
an interval $E \subset\mathbb{R}$.
It means that $X$ is a regular Feller process having jumps only in one
direction which we assume, without loss of generality,
to be of the spectrally negative type. That is, $X$ does not have
positive jumps;
$\mathbb{P}_x(\sup_{t\geq0} (X_t - X_{t-}) > 0 ) = 0, \forall x\in E $.
For sake of simplicity, we also assume that the process $X$ has
infinite lifetime.
We wish to extend the results of the previous section to this class of
stochastic processes.
In comparison to the diffusion case, the difficulty is that we do not
have, in general,
any information\vadjust{\goodbreak} on the excessive functions for this class of Markov
processes. Indeed,
the infinitesimal generator associated to $X$ is an
integro-differential linear operator
for which there does not exist general results regarding the solutions
to the boundary
value problem (\ref{Sturmliouville}). Nevertheless, as explained in the
following,
the one-sided feature of $X$ allows us to identify increasing excessive
functions.
\begin{proposition} \label{propsn}
For any $q>0$, there exists a unique increasing left-continuous
function $h^+_q \dvtx E \mapsto[0,\infty] $, such that, for any $x, y \in
E$ with $x \leq y$,
\[
\mathbb{E}_x [e^{-q T_y}] = \frac{h^+_q(x)}{h^+_q(y)}.
\]
In the case $X$ is recurrent, the function $ x \mapsto h^+_q(x)$ is continuous.
\end{proposition}
\begin{pf}
As a consequence of the regularity assumption, it is well known (see, e.g.,
\cite{Blumenthal-Getoor-68})
that for
each singleton $\{y\} \in E$, X admits a local time at $y$, which we
denote by
$L^y=(L^y_t)_{t\geq0}$. The continuous additive functional $L^y$ is
determined by
its $q$-potential, $u_q^y$, which is finite for any $q>0$ and given by
\[
u_q^y(x)=\mathbb{E}_x \biggl[\int_0^{\infty} e^{-q t} \,dL_t^y
\biggr].
\]
From the definition of $L^y$ and the strong Markov property, we
obtain the identity (see~\cite{Blumenthal-Getoor-68}, Chapter V.3)
\[
\mathbb{E}_x [e^{-q T_y}] = \frac{u_q^y(x)}{u_q^y(y)},\qquad
x,y \in E.
\]
Next, following It\^o and McKean
\cite{ItoMckean} or Pitman and Yor~\cite{Pitman-Yor-81}, for instance,
we write $\phi_q(x,y) = \mathbb{E}_x [e^{-q T_y}]$ and for a fixed
$z_0\in E$, we define
\[
h_q^+(y) = \cases{
\phi_q(y,z_0), &\quad $y\leq z_0$,\vspace*{2pt}\cr
1/\phi_q(z_0,y), &\quad $y> z_0$.}
\]
Next, using the fact that $X$ has no positive jumps, we get for any
$x\leq z\leq y$, and by means of the
strong Markov property,
\[
\phi_q(x,y) = \phi_q(x,z)\phi_q(z,y).
\]
Thus, from the identity
\[
\phi_q(x,z) = h_q^+(x)/h_q^+(z)
\]
we deduce that the choice of the reference point $z_0$ affects $h_q^+$
only by a~constant factor.
The monotonicity of the mapping $h_q^+$ follows readily from its
definition and the absence
of positive jumps for $X$.
We recall that $x\mapsto\mathbb{E}_x [e^{-q T_y}] $ is a
$q$-excessive function
(see, e.g.,~\cite{Blumenthal-Getoor-68}, page 74).
By linearity, the mapping $x\mapsto h^+_q(x)$ is also $q$-excessive
and thus finely continuous. So, the Feller property of $X$
implies that the increasing excessive function $h^+_q$ is lower
semi-continuous (see
\cite{hawkes}) and hence, left-continuous.
Then, the claim of the last assertion is a straightforward
consequence of the fact that if $X$\vadjust{\goodbreak} is also
recurrent then the fine
topology coincides with the initial topology of~$E$ (see, e.g.,
\cite{Bally-Stoica-87}, page 243).
\end{pf}
We point out that, in Patie and Vigon~\cite{patievigon},
Proposition~\ref{propsn} is extended to a~larger class of homogenous
Markov processes with only negative jumps.
Following the proof of Theorem~\ref{thmmain1}, and observing that
$h^{(-\infty, u)}(u) = h_q^+(u)/\break h_q^+(x)$,
the proof of the theorem below goes through verbatim.
\begin{theo} \label{thmmain2}
With the notation used in Theorem~\ref{thmmain1} we assume that
\[
D^{*}:=\sup_{u \in E} \frac{ g(u)+\delta(u)}{ h^+_q(u)} =
\sup_{u\geq x, u \in E} \frac{ g(u)+\delta(u)}{ h^+_q(u)} <\infty.
\]
Then, for any $x \in E$
\[
\sup_{T \in\Sigma_\infty^X}
\mathbb{E}_x[e^{-qT}g(X_T)-C_T] =D^*h^+_q(x)-\delta(x).
\]
If there exists a point $ u^{*} \geq x$ such that
$D^{*}= \frac{ g(u^*)+\delta(u^*)}{ h^+_q(u^*)}$, then the
optimal stopping time of the problem (\ref{eqoptstop}) is given by
\[
T_{u^*} = \inf\{t >0; X_t \geq u^* \}.
\]
\end{theo}
\begin{remark}
Since Lemma~\ref{timechangelema} applies also in this more general
framework, one could have also considered a
random discounting factor in the previous result.
\end{remark}

\section{Examples}
We now illustrate our methodology by presenting the solutions to some
new optimal stopping problems. We consider
both the diffusion case and also the case when the processes are of the
spectrally negative type.

\subsection{An optimal stopping problem with cost of observations for one-sided
L\'evy processes}
Let $Z=(Z_t)_{t\geq0}$ be a spectrally negative L\'evy process
starting from $x \in\mathbb{R}$, that is, a
process with stationary and independent increments having only negative jumps.
Plainly, the law of $Z$ is characterized by~$\psi$, the Laplace
exponent of $Z_1$, which admits the following
L\'evy--Khintchine representation, for any $u\geq0$,
%
\begin{equation} \label{eqlap-lev}
\psi(u)= \frac{1}{2} \sigma^2 u^2 + b u + \int_{-\infty}^0 \bigl( e ^{ux}-
1 - ux \mathbh{1}_{\{-1 < x \leq0\}}\bigr) \nu(dx),
\end{equation}
where $b\in\mathbb{R}$ and $\sigma\geq0$ and the measure $\nu$ is
such that $ \int_{-\infty}^0
(1\wedge y^2) \nu(dy) < \infty$. Next, recalling that $\psi$ is a
convex function on $[0,\infty)$
with $\lim_{u\rightarrow\infty} \psi(u)=+\infty$, we denote by
$\theta
$ the nonnegative largest
root of the equation $\psi(u)=0$. We also mention that
being continuous and increasing on $[\theta,\infty)$, $\psi$ has
a~well-defined inverse function $\phi\dvtx[0,\infty)\rightarrow
[\theta,\infty)$ which is also continuous and increasing. We refer to
the excellent monographs of
Bertoin~\cite{Bertoin-96} and Kyprianou~\cite{Kyprianou-06} for
background on L\'evy processes.
We now consider a perpetual American option\vadjust{\goodbreak} in a market driven by $Z$
where the
agent takes into account some costs from hedging, which might come
from some transaction costs or liquidity issues. We assume that
the agent hedges continuously and that he chooses the cost of observations
$c(x) = \exp(\gamma x)$ to be of exponential form.
The payoff at time $t> 0$ of such a product can be written as
%
\begin{equation}
e^{-q t}g(e^{Z_t}) - \int_0^{t} e^{\gamma Z_s} e^{-qs} \,ds,
\end{equation}
where $q>0$ is the risk-free rate, $\gamma> -1$ and $g$ is a smooth
function. Next, assuming that $q=\psi(1)$
it is easy to check that one may choose the probability $\mathbb{P}_x$
as the risk neutral probability measure.
Therefore, in the sequel, we suppose that the characteristics $b,\sigma
$ and $\nu$ are chosen such that
%
\begin{equation}\label{martingmea}
q=\psi(1).
\end{equation}

\subsubsection{The Brownian motion with drift case}
We start with the case where $Z_t= W_t^{(b)}= b t +W_t$ is a Brownian
motion with drift $b$ starting from~$0$ and
the reward function $g$ is defined, for any $0 \leq L \leq K$, by
%
\begin{equation}
g(x)=
\cases{
L - x, &\quad if $x\leq L$,\cr
0, &\quad if $L \leq x \leq K$,\cr
x -K, &\quad if $x \geq K$.}
\end{equation}
This function corresponds to the payoff function of a strangle option
which is a combination of a put with
exercise price $L$ and a call with exercise price $K$. In this case,
the condition (\ref{martingmea})
is fulfilled if $q=\frac{1}{2} +b$.
We want to compute the constant
\[
V_W = \sup_{T \in\Sigma_ \infty} \mathbb
{E}_0\biggl[e^{-qT}g\bigl(e^{W_T^{(b)}}\bigr)-\int_0^T
e^{\gamma W_s^{(b)}}e^{-qs}\,ds \biggr].
\]
The case without cost of observations, that is, $ c\equiv0$, has
already been studied by Beibel and Lerche
\cite{Beibel-Lerche-97}.
Next, it is well known that the functions $h^+_q$ and~$h_q^-$ defined
in (\ref{defnitfi}) are
given by
\[
h^+_q(x)=D_1 e^{\alpha_1 x},\qquad h_q^-(x)= D_2 e^{\alpha_2 x},
\]
where $
\alpha_1= -b + \sqrt{2q+b^2}$ and $ \alpha_2= -b - \sqrt{2q+b^2} $,
and $D_1$, $D_2$ are positive real numbers.\vadjust{\goodbreak} The function $\delta$,
defined in (\ref{potentiel}), is finite if $ q > b \gamma+
\frac{\gamma^2}{2}$ and is given by
\[
\delta(x)= \frac{e^{\gamma x}}{q-\gamma b - \gamma^2 /2},\qquad x\in
\mathbb{R}.
\]
Next, we introduce the function
\[
G_p(x)=\frac{g(x)+\delta(x)}{p e^{\alpha_1 x} + (1-p) e^{\alpha_2 x}},
\]
which, according to Theorem~\ref{thmmain1}, gives the solution to our problem.
Then, if $q>\gamma b +\frac{\gamma^2 }{2}$ and $ b > -1$ then
$\alpha
_1 >1$ and $\alpha_2 < -1$. Thus, we\vadjust{\goodbreak}
verify easily the following inequalities:
%
\begin{equation}\label{BL2}
\sup_{ x \leq0}\bigl( e^{-\alpha_1 x} \bigl(\delta(x)+g(x)\bigr)\bigr)
 > \sup_{ x \geq
0}\bigl( e^{-\alpha_1 x}
\bigl(\delta(x)+g(x)\bigr)\bigr) >0
\end{equation}
and
%
\begin{equation}\label{BL3}
\sup_{ x \geq0}\bigl( e^{-\alpha_2 x} \bigl(\delta(x)+g(x)\bigr)\bigr) > \sup_{ x \leq
0}\bigl( e^{-\alpha_2 x}
\bigl(\delta(x)+g(x)\bigr)\bigr) >0.
\end{equation}
Note that
\[
\sup_{x\geq0} G_p(x)= \sup_{x\geq0; \delta(x)+g(x)>0} G_p(x)= \biggl(
\inf_{x\geq0;
\delta(x)+g(x)>0} \frac{p e^{\alpha_1 x} + (1-p) e^{\alpha_2
x}}{g(x)+\delta(x)}\biggr)^{-1}
\]
and
\[
\sup_{x\leq0} G_p(x)= \sup_{x\leq0; \delta(x)+g(x)>0} G_p(x)= \biggl(
\inf_{x\leq0;
\delta(x)+g(x)>0} \frac{p e^{\alpha_1 x} + (1-p) e^{\alpha_2
x}}{g(x)+\delta(x)}\biggr)^{-1}.
\]
Then for all $p \in(0,1)$
\[
0< \sup_{x\geq0; \delta(x)+g(x)>0} G_p(x)\leq\frac{1}{p} \sup
_{x\geq0} \bigl(e^{-\alpha_1 x}
\bigl(g(x)+\delta(x)\bigr) \bigr) < +\infty
\]
and
\[
0< \sup_{x\leq0; \delta(x)+g(x)>0} G_p(x)\leq\frac{1}{1-p} \sup
_{x\leq0} \bigl(e^{-\alpha_2 x}
\bigl(g(x)+\delta(x)\bigr) \bigr) < +\infty.
\]
We assume that $
\log(L) \leq0= W^{b}_0 \leq\log(K), q>1/2 +b$ and $ b >-1$.
Then~(\ref{BL2}) and (\ref{BL3}) hold and, as in Beibel and Lerche
\cite{Beibel-Lerche-97}, Lemma 1, page 98,
there exists a number $p^* \in(0,1)$ such that
\[
\sup_{x\geq0} G_{p^*}(x)= \sup_{x\leq0}G_{p^*}(x).
\]
Let $ x_1, x_2$ and $p^*$ be solutions with $ x_1 >\log K, x_2 < \log
L $ and $p^* \in(0,1)$ of
the following system:
%
\begin{equation}
\cases{
\displaystyle \frac{(q-\gamma b - \gamma^2 /2)(e^{x_1}-K )+e^{-x_1}}{p^* e^{\alpha_1
x_1} + (1-p^*) e^{\alpha_2 x_1}}=
\frac{(q-\gamma b - \gamma^2 /2)(L- e^{x_2} )+e^{-x_2}}{p^*
e^{\alpha
_1 x_2} + (1-p^*) e^{\alpha_2
x_2}},
\vspace*{2pt}\cr
\displaystyle \frac{e^{x_2}(\gamma b + \gamma^2 /2-q)- e^{-x_2}}{(q-\gamma b -
\gamma
^2 /2)(L-e^{x_2})+e^{-x_2}}=
\frac{p^* \alpha_1 e^{\alpha_1 x_2} + (1-p^*)\alpha_2 e^{\alpha_2
x_2}}{p^* e^{\alpha_1 x_2} + (1-p^*)
e^{\alpha_2 x_2}},\vspace*{2pt}\cr
\displaystyle \frac{e^{x_1}(q-\gamma b - \gamma^2 /2)
-e^{-x_1}}{(e^{x_1}-K)(q-\gamma b - \gamma^2 /2) +e^{-x_1}} =
\frac{p^* \alpha_1 e^{\alpha_1 x_1} + (1-p^*)\alpha_2 e^{\alpha_2
x_1}}{p^* e^{\alpha_1 x_1} + (1-p^*) e^{\alpha_2
x_1}}.}\hspace*{-29pt}
\end{equation}
Let
\[
M^*= \frac{(q-\gamma b - \gamma^2 /2)(e^{x_1}-K )+e^{-x_1}}{p^*
e^{\alpha_1 x_1} + (1-p^*) e^{\alpha_2 x_1}}
= \frac{(q-\gamma b - \gamma^2 /2)(L- e^{x_2} )+e^{x_2}}{p^*
e^{\alpha
_1 x_2} + (1-p^*) e^{\alpha_2 x_2}}.
\]
Then,
\[
\sup_{T\in\Sigma^X_\infty}\mathbb{E}_0\biggl[e^{-qT}g\bigl(W_T^{(b)}\bigr)-\int_0^T
c\bigl(W_s^{(b)}\bigr)e^{-qs}\,ds \biggr]= M^* - \frac{1}{q-\gamma b - \gamma^2 /2}\vadjust{\goodbreak}
\]
and the optimal stopping time is
\[
T_{(x_1,x_2)}= \inf\bigl\{ t>0; W^{(b)}_t \notin(x_1,x_2)\bigr\}.
\]

\subsubsection{The spectrally negative L\'evy case}
In the second example, we consider as the reward function
\[
g(x)= (x-K)^+.
\]
We write, for any $x\in E$,
%
\begin{equation}\label{Levyexple}
\mathcal{V}_Z(x)=\sup_{T \in\Sigma_ \infty} \mathbb{E}_x \biggl[ e ^{-qT}
( e^{Z_T}-K )^+ - \alpha
\int_0^T e^{\gamma Z_s} e^{-qs} \,ds \biggr],
\end{equation}
where $\alpha>0$ and we recall that we choose $q=\psi(1)$.
\begin{proposition}
We assume that $p_{\gamma} = \psi(1) - \psi(\gamma)>0$ and
%
\begin{equation} \label{eqcx}
x < x^*=\frac{1}{\gamma} \log\biggl(\frac{p_{\gamma}K}{(1-\gamma)\alpha}\biggr).
\end{equation}
Then
\[
\mathcal{V}_Z(x)= e^{x-x^*} \biggl( (e^{x^*}-K)^+ +
\alpha\frac{ e ^{\gamma x^*}}{p_{\gamma} } \biggr)-\alpha\frac{ e
^{\gamma x}}{p_{\gamma} }
\]
and the optimal stopping time is $T_{x^*}= \inf\{t>0; Z_t\geq x^*\}$.
\end{proposition}
\begin{pf}
First, by means of Fubini's theorem and using the fact that $p_{\gamma
}>0$, we easily get that, for any $ x\in
E$,
\[
\delta(x) = \alpha\mathbb{E}_x \int_0 ^\infty e^{-qt} e^{\gamma Z_t}
\,dt = \alpha e^{\gamma x} \int_0 ^\infty e^{-qt+ t \psi(\gamma)} \,dt
= \alpha\frac{ e ^{\gamma x}}{p_{\gamma} } .
\]
Then, from Lemma~\ref{ramenbabel}, we deduce that
\[
V_Z(x)= \sup_{T \in\Sigma_ \infty} \mathbb{E}_x \bigl[ e ^{-qT} \bigl(
(e^{Z_T}- K )^+
+ \delta(Z_T) \bigr) \bigr] - \delta(x).
\]
Next, we recall from Bertoin
\cite{Bertoin-96}, pages 189 and 190, that, for any $x\leq y$,
\[
\mathbb{E}_x [ e^{-r T_y} ]= e^{- \phi(r) (y-x)},\qquad r\geq0.
\]
Then, writing
\[
G(u)= \frac{ ( e^{u}-K )^+ + \delta(u) }{e^{u-x}},
\]
we have that
\[
G'(u)= e^{x-u}\biggl(K \mathbh{1}_{\{u \geq\log K\}} + \frac{\alpha
(\gamma
- 1)}{p_{\gamma}} e^{\gamma u}\biggr).
\]
Since $G'(u) \geq0$ on $(-\infty,x^*]$ and $G'(u) \leq0$ otherwise,
we deduce from (\ref{eqcx}) that
\[
\sup_{u \in E} G(u) =\sup_{u \geq x }G(u) < \infty.
\]
The proof of the claim follows then by applying Theorem~\ref{thmmain2}.
\end{pf}

\subsection{Optimal stopping problems associated to self-similar
positive Markov processes of the spectrally negative type}
Let $Z=(Z_t)_{t\geq0}$ again be a spectrally negative L\'evy process.
We introduce,
for any $\alpha>0$, the process $X=(X_t)_{t\geq0}$ defined by
%
\begin{equation} \label{eqlamptransf}
\log(X_t) = Z_{A_t},\qquad t\geq0,
\end{equation}
where
\[
A_t = \inf
\biggl\{ s \geq0; \Sigma_s = \int_0^s e^{\alpha Z_r} \,dr
> t \biggr\}.
\]
We denote its law by $\mathbb{Q}_{x}$ when it starts from $x>0$.
Lamperti~\cite{Lamperti-72} showed that $X$ is an $\alpha
$-self-similar positive
Markov process on $(0,\infty)$,
that is, a~Feller process which enjoys the
following $\alpha$-self-similarity property, for any $c>0$,
%
\begin{equation} \label{eqself}
((X_{c^{\alpha}t})_{t\geq0}, \mathbb{Q}_{x})
\stackrel{(d)}{=}((cX_{t})_{t\geq0},
\mathbb{Q}_{x}).
\end{equation}
It is plain that $X$ is also of the spectrally
negative type, in the sense that it has no positive jumps. Next, we
recall that the law of $Z$ is characterized by
its Laplace exponent, $\psi$, which is of the form (\ref{eqlap-lev}).
In the sequel, writing $\theta$ for the largest root in $[0,\infty)$ of
the equation $\psi(u) = 0$,
we assume that the following conditions
%
\begin{equation} \label{eqcond}
\theta<\alpha\quad\mbox{and}\quad \lim_{u\rightarrow\infty} \frac{\psi
(u)}{u} = \infty
\end{equation}
hold. The first condition secures that the lifetime of $X$ is infinite
since in the case $0<\theta<\alpha$, we consider $X$ to be the unique
recurrent extension which hits and leaves $0$ continuously (see Rivero
\cite{Rivero-05} for more details). Under the second condition, the
paths of the process $X$ are of unbounded variation on any compact
interval and the process $X$ is regular. Next, we introduce more
notation taken from Patie~\cite{Patie-06c}. Define for any integers $n$
\[
\ant^{-1}=\prod_{k=1}^n \psi(\alpha k) ,\qquad a_0(\psi,\alpha)=1,
\]
and we introduce the entire function $\Ip$ which admits the series
representation
\[
\Ip_{\psi,\alpha}(z)=\sum_{n=0}^{\infty} \ant z^{n},\qquad
z \in\mathbb{C}.
\]
It is important to note that whenever $\theta$, the largest root of the
equation $\psi(u)=0$,
satisfies $\theta<\alpha$, it follows that all of the coefficients in
the definition of $\Ip_{\psi,\alpha}$\vadjust{\goodbreak}
are strictly positive. Then, Patie~\cite{Patie-06c}, Theorem 2.1,
characterized the Laplace transform of
\[
T_a=\inf\{t>0; X_t \geq a\}
\]
as follows.
Suppose that $0 \leq x \leq a$. Then, for any $q\geq0$,
we have
%
\begin{equation} \label{eqlt}
{\mathbb{E}}_x [e^{-q T_a } ]
= \frac{\Ip_{\psi,\alpha}(qx^{\alpha})}{\Ip_{\psi,\alpha
}(qa^{\alpha})}.
\end{equation}
Next, we introduce the Ornstein--Uhlenbeck process associated to $X$
which is defined, for any $t\geq0$, by
%
\begin{equation}\label{eqdefp}
U_t =e'_{\lambda}(-t)X_{e_{\chi}(t)},\qquad t \geq0,
\end{equation}
where $e_{\lambda}(t)= \frac{e^{\lambda t}-1}{\lambda}$,
$\chi=\alpha\lambda$ and we write $v_{\lambda}(t)= \frac{\log
(1+\lambda t)}{\lambda}$ the continuous
increasing inverse function of $e_{\lambda}$. These processes were
introduced and studied by Carmona, Petit and Yor
\cite{Carmona-Petit-Yor-97}. In particular, they proved that $U$ is
a~Feller process and under the
conditions (\ref{eqcond}), $U$ has also infinite lifetime and is regular.
In~\cite{Patie-ouq-06}, the author computed the Laplace transform of the
first passage times above of $U$ as follows.
With the obvious notation, for any $r \geq0$ and $0\leq x\leq a$,
%
\begin{equation}
\label{eqlapssm}
{\mathbb{E}}_x [e^{-r T^U_a } ] =
\frac{\mathcal{I}_{\psi,\alpha}(r;x^{\alpha})}{\mathcal{I}_{\psi
,\alpha
}(r;a^{\alpha})},
\end{equation}
where
%
\begin{equation} \label{eqf1}
\mathcal{I}_{\psi,\alpha}(q;x)= \sum_{n=0}^{\infty} \frac{\Gamma
(q+n)}{\Gamma(q)} a_n(\psi;\alpha) x^{n}
\end{equation}
and $\Gamma$ stands for the gamma function. By means of classical
criteria on power series, it is easily seen that,
under the second condition in (\ref{eqcond}),
that $\mathcal{I}_{\psi,\alpha}(q;x)$ is an entire function in $x$ and
is analytic on
the domain ${\{q\in\mathbb{C}; \Re(q)>-1\}}$. We shall also need the
following representation of the function $\mathcal{I}_{\psi,\alpha}(q;x)$
%
\begin{equation} \label{eqf2}
\mathcal{I}_{\psi,\alpha}(q;x)= \frac{1}{\Gamma(q)}\int_0^{\infty
}\mathcal{I}_{\psi,\alpha}(rx)e^{-r}r^{q-1}\,dr,
\end{equation}
which is readily obtained by using the integral representation of the
gamma function $\Gamma(q)=\int_0^{\infty} e^{-r}r^{q-1}\,dr,
\Re(q)>0$. We postpone to the end of this section the description of
some specific examples of these power series. Let us assume that $g$
is a continuous function and $q,\beta>0$. We are ready to introduce the
following optimal stopping problems:
\begin{eqnarray*}
\mathcal{V}^{X}_{g}(x) &\stackrel{\Delta}{=}& \sup_{T \in\Sigma
_\infty
^X} \mathbb{E}_x[e^{-qT}g(X_T)] ,\\
\mathcal{V}^{U}_{g}(x) &\stackrel{\Delta}{=}& \sup_{T \in\Sigma
_{\infty
}} \mathbb{E}_x[e^{-qT}g(U_T)],\\
\mathcal{V}^{U,\bigtriangleup}_g(x) &\stackrel{\Delta}{=}& \sup_{T
\in
\Sigma_{\infty}} \mathbb{E}_x[e^{-q\bigtriangleup_T}
g(U_T)],\qquad \bigtriangleup_{t} =\int_0^t U_s^{-\alpha} \,ds,\\
\mathcal{V}^{S,q}_{g}(x) &\stackrel{\Delta}{=}& \sup_{T \in\Sigma
_{\infty}} \mathbb{E}_x\biggl[e^{-qT}g\biggl(\frac{e^{\alpha Z_T}}{1+\beta
\int_0^Te^{\alpha Z_s}\,ds}\biggr)\biggr].
\end{eqnarray*}
We note that in the case $Z$ is a Brownian motion with drift, the last
optimal stopping problem is intimately connected to the so-called
integral option problem studied by Kramkov and Mordetski
\cite{Kramkov-Morde-94}.
\begin{proposition}
Let us write for some function $h$
\[
a_*(h(\cdot)) = \arg\max_{u\geq x, u \in E} \frac{g(u)}{h(u^{\alpha})}.
\]

\begin{longlist}[(1)]
\item[(1)] If $a_*=a_* (\mathcal{I}_{\psi,\alpha}(q\cdot))$ exists and $x<a_*
$ then
\[
\mathcal{V}^{X}_g(x) = \frac{\mathcal{I}_{\psi,\alpha}(q x^{\alpha
})}{\mathcal{I}_{\psi,\alpha}(qa_*^{\alpha})} g(a_*).
\]
\item[(2)] If $a_*=a_*(\mathcal{I}_{\psi,\alpha}(\frac{q}{\chi};\chi\cdot))$
exists and $x<a_* $ then
\[
\mathcal{V}^{U}_g(x) = \frac{\mathcal{I}_{\psi,\alpha}
({q}/{\alpha};\chi x^{\alpha})}{\mathcal{I}_{\psi,\alpha
}({q}/{\alpha};\chi a_*^{\alpha})}
g(a_*).
\]
\item[(3)] If $a_*=a_*(\mathcal{I}_{\psi,\alpha}(\frac{q}{\chi};\chi\cdot))$
exists and $x<a_* $ then
\[
\mathcal{V}^{U,\bigtriangleup}_g(x) = \biggl(\frac{x}{a_*}\biggr)
^{\gamma}\frac
{\mathcal{I}_{\psi,\alpha}({\gamma}/{\alpha};\chi x^{\alpha
})}{\mathcal{I}_{\psi,\alpha}({\gamma}/{\alpha};\chi
a_*^{\alpha})} g(a_*).
\]
\item[(4)] If $a_*=a_*(\mathcal{I}_{\psi,\alpha}(\frac{q}{\chi};\chi\cdot))$
exists and $x<a_* $ then
\begin{eqnarray*}
\mathcal{V}^{S,q}_g(x) & = & \Psi^U_{\bigtriangleup,g}(1)\\
&=& \frac{\mathcal{I}_{\psi,\alpha}({q}/{\chi};\chi x^{\alpha})}
{\mathcal{I}_{\psi,\alpha}({q}/{\chi};\chi a_*^{\alpha})} g(a_*).
\end{eqnarray*}
\end{longlist}
In all the above cases the optimal stopping time is given by $T_{a_*}$.
\end{proposition}
\begin{pf}
The first item follows readily from the identity (\ref{eqlapssm}) and
Theorem~\ref{thmmain2}. Next, let
\[
{T^X_{y,\alpha} =
\inf\{t>0; X_t = y(1+\alpha\lambda t)^{{1}/{\alpha}}\}}.
\]
The Mellin transform of the positive random variable $T^X_{y,\alpha}$
has been computed by Patie (see~\cite{Patie-ouq-06}, Theorem 2).
However, for sake of completeness, we provide a
slightly different proof here which relies on a device introduced by Shepp
\cite{Shepp-67}.
In the proof of~\cite{Patie-06c}, Theorem 1, it is shown that the
mapping $x\mapsto\mathcal{I}_{\psi,\alpha}(rx^{\alpha})$
is an $r$-eigenfunction for the infinitesimal generator of~$X$.
Thus, by the Dynkin formula, using the fact that the function $\mathcal
{I}_{\psi,\alpha}$
is increasing\vadjust{\goodbreak} and applying the dominated convergence theorem, we
deduce that
\[
{\mathbb{E}}_x \bigl[e^{-r T^X_{y,\alpha}}\mathcal{I}_{\psi,\alpha
}\bigl(ry^{\alpha}(1+\alpha\lambda T^X_{y,\alpha})\bigr) \bigr] =\mathcal{I}_{\psi
,\alpha}(rx^{\alpha}).
\]
Integrating both sides of the previous identity
by the measure $e^{-\chi r}
r^{{q}/{\chi}-1} \,dr$, using Fubini's theorem and the change of
variable $u=r( T^X_{y,\alpha} + \chi)$, we get
\[
\mathbb{E}_x [(1+\chi T^X_{y,\alpha} )^{-{q}/{\chi}} ] =\frac
{\mathcal{I}_{\psi,\alpha}({q}/{\chi};\chi x^{\alpha})}
{\mathcal{I}_{\psi,\alpha}({q}/{\chi};\chi a^{\alpha})},
\]
where we have used (\ref{eqf2}). Moreover, from the definition of $U$
(\ref{eqdefp}), we observe the identity
%
\begin{equation} \label{eqist}
T^U_y = v_{\chi}(T^X_{y,\alpha}) \qquad\mbox{a.s.}
\end{equation}
Hence, we obtain that
\[
\mathbb{E}_x [e^{-q T_a^{U} } ] =\frac{\mathcal{I}_{\psi,\alpha
}({q}/{\chi};\chi x^{\alpha})}
{\mathcal{I}_{\psi,\alpha}({q}/{\chi};\chi a^{\alpha})}.
\]
We deduce the item (2) by an application of Theorem~\ref{thmmain2}.
We complete the proof of the proposition from~\cite{Patie-ouq-06},
Corollary 3.2.
\end{pf}
\begin{remark}
We note that the identity (\ref{eqist}) combined with the item (2) of
the previous proposition
allows us to solve the following nonhomogeneous optimal stopping time problem
\[
\mathcal{V}^{X}_{g, \alpha}(x) = \sup_{T \in\Sigma_\infty^X}
\mathbb
{E}_x\biggl[(1+2\lambda T)
^{-q}g\biggl(\frac{X_T}{(1+2 \lambda T)^{{1}/{\alpha}}}\biggr)\biggr].
\]
Indeed, we easily deduce that
\[
\mathcal{V}^{X}_{g,\alpha}(x) = \frac{\mathcal{I}_{\psi,\alpha
}({q}/{\chi};\chi x^{\alpha})}
{\mathcal{I}_{\psi,\alpha}({q}/{\chi};\chi a_*^{\alpha})} g(a_*),
\]
where $a_*$ is characterized by
\[
a_* = \arg\max_{u\geq x, u \in E} \frac{g(u)}{\mathcal{I}_{\psi
,\alpha
}({q}/{\chi};\chi u^{\alpha})}.
\]
\end{remark}

In what follows, we provide some examples of the power series and we
refer to~\cite{Patie-06c} for the description of additional examples.

\subsubsection*{The modified Bessel functions}
We consider $Z$ to be a Brownian motion with drift $\nu\geq0$, that is,
$\psi(u)=\frac{1}{2}u^2+\nu u$ and we set $\alpha=2$. Its associated
self-similar process is well known to be a Bessel process of index
$\nu$. We have
\[
a_n(\psi;2)^{-1}
= 2^n n! \frac{\Gamma( n- \nu+1)}{\Gamma(-\nu+1)},\qquad a_0 =1.
\]
Thus, we get
\[
\mathcal{I}_{2,\psi}(x)=(x/2)^{\nu/2}\Gamma(-\nu+1){\mathrm{I}}_{-\nu
}\bigl(\sqrt{2x}\bigr),\vadjust{\goodbreak}
\]
where\vspace*{1pt} ${\mathrm{I}}_{\nu}(x) = \sum_{n=0}^{\infty}\frac{(x/2)^{\nu
+2n}}{n!\Gamma(\nu+n+1)}$
stands for the modified Bessel function of index~$\nu$ (see, e.g.,
\cite{Lebedev-72}, Chapter 5) and
\[
\mathcal{I}_{2,\psi}(q;x^2) = \Phi\biggl(q,1-\nu,\frac{x^2}{2}\biggr),
\]
where $\Phi(q,\nu,x)=\sum_{n=0}^{\infty}\frac{(q)_n}{(\nu)_n n!}x^n$
stands for the confluent hypergeometric function of the first kind
(see, e.g.,
\cite{Lebedev-72}, Section 9.9, page 260)
and $(q)_{n}=\frac{\Gamma(q+n)}{\Gamma(n)}$ stands for the
Pochhammer symbol.\vspace*{-3pt}

\subsubsection*{Some generalized Mittag--Leffler functions} \label{secml}
In~\cite{Patie-06-poch}, the author introduced a parametric
family of one-sided L\'evy processes which are characterized by the
following Laplace exponent, for any $ 1<\alpha< 2$,
and $\gamma
>1-\alpha$,
%
\begin{equation} \label{eqlappoch}
\psi_{\gamma}(u)= \bigl((
u+\gamma-1)_{\alpha}-(\gamma-1)_{\alpha}\bigr).
\end{equation}
Its L\'evy measure is absolutely continuous with a density $f$ given by
\[
f(y)=C\frac{e^{(\alpha+\gamma-1)y}}{(1-e^{y})^{\alpha+1}},\qquad
y<0,
\]
where $C$ is a positive constant.
We focus on the case $\gamma=1$ in (\ref{eqlappoch}). We have $ \psi
_{1}(u)=\psi(u)=
(u)_{\alpha}$ and
\[
a_n(\psi;\alpha)^{-1}=\frac{\Gamma(\alpha(n+1))}{\Gamma(\alpha)},\qquad
a_0=1.
\]
Thus, the power series can be written as
\begin{eqnarray*}
\mathcal{I}_{\psi,\alpha}(x)&=&\Gamma(\alpha)\mathcal{M}_{\alpha
,\alpha
}(\alpha
x) ,\\
\mathcal{I}_{\psi,\alpha}(q;x)&=&\Gamma(\alpha)\mathcal
{M}^q_{\alpha
,\alpha}(\alpha
x),
\end{eqnarray*}
where $ \mathcal{M}_{\alpha,\beta}(x)= \sum_{n=0}^{\infty} \frac{
x^n}{\Gamma(\alpha n+\beta)}$ [resp., $ \mathcal{M}^q_{\alpha,\beta
}(x)= \sum_{n=0}^{\infty} \frac{(q)_n x^n}{\Gamma(\alpha n+\beta)}$]
stands for the Mittag--Leffler function of parameters
$\alpha,\beta>0$ (resp., of parameters $\alpha,\beta,q>0$) introduced
by Prabhakar~\cite{Prabhakar-71}.\vspace*{-3pt}

\section*{Acknowledgments}
The authors are indebted to an anonymous referee for the insightful
comments that helped in improving the presentation of the paper. This
research has been carried out within the NCCR FINRISK project on
``Credit Risk and Non-Standard Sources of Risk in Finance.'' Financial
support by the National Center of Competence in Research Financial
valuation and Risk Management (NCCR FINRISK) is gratefully
acknowledged. NCCR-FINRISK is a research program supported by the Swiss
National Science Foundation.\vspace*{-3pt}


%

\printaddresses


\begin{thebibliography}{36}

\bibitem{Alili-Kyprianou-05}
\begin{barticle}[mr]
\bauthor{\bsnm{Alili},~\bfnm{L.}\binits{L.}} \AND
  \bauthor{\bsnm{Kyprianou},~\bfnm{A.~E.}\binits{A.~E.}}
(\byear{2005}).
\btitle{Some remarks on first passage of {L}\'evy processes, the {A}merican put
  and pasting principles}.
\bjournal{Ann. Appl. Probab.}
\bvolume{15}
\bpages{2062--2080}.
\bid{doi={10.1214/105051605000000377}, issn={1050-5164}, mr={2152253}}
\bptok{imsref}%
\end{barticle}\vadjust{\goodbreak}
\endbibitem

\bibitem{Bally-Stoica-87}
\begin{barticle}[mr]
\bauthor{\bsnm{Bally},~\bfnm{Vlad}\binits{V.}} \AND
  \bauthor{\bsnm{Stoica},~\bfnm{Lucretiu}\binits{L.}}
(\byear{1987}).
\btitle{A class of {M}arkov processes which admit local times}.
\bjournal{Ann. Probab.}
\bvolume{15}
\bpages{241--262}.
\bid{issn={0091-1798}, mr={0877600}}
\bptok{imsref}%
\end{barticle}
\endbibitem

\bibitem{Baurdoux-07}
\begin{barticle}[mr]
\bauthor{\bsnm{Baurdoux},~\bfnm{E.~J.}\binits{E.~J.}}
(\byear{2007}).
\btitle{Examples of optimal stopping via measure transformation for processes
  with one-sided jumps}.
\bjournal{Stochastics}
\bvolume{79}
\bpages{303--307}.
\bid{doi={10.1080/17442500600856297}, issn={1744-2508}, mr={2308078}}
\bptok{imsref}%
\end{barticle}
\endbibitem

\bibitem{Beibel-Lerche-97}
\begin{barticle}[mr]
\bauthor{\bsnm{Beibel},~\bfnm{M.}\binits{M.}} \AND
  \bauthor{\bsnm{Lerche},~\bfnm{H.~R.}\binits{H.~R.}}
(\byear{1997}).
\btitle{A new look at optimal stopping problems related to mathematical
  finance}.
\bjournal{Statist. Sinica}
\bvolume{7}
\bpages{93--108}.
\bid{issn={1017-0405}, mr={1441146}}
\bptok{imsref}%
\end{barticle}
\endbibitem

\bibitem{Beibel-Lerche-00}
\begin{barticle}[mr]
\bauthor{\bsnm{Beibel},~\bfnm{M.}\binits{M.}} \AND
  \bauthor{\bsnm{Lerche},~\bfnm{H.~R.}\binits{H.~R.}}
(\byear{2000}).
\btitle{A note on optimal stopping of regular diffusions under random
  discounting}.
\bjournal{Teor. Veroyatn. Primen.}
\bvolume{45}
\bpages{657--669}.
\bid{doi={10.1137/S0040585X9797852X}, issn={0040-361X}, mr={1968720}}
\bptok{imsref}%
\end{barticle}
\endbibitem

\bibitem{Bertoin-96}
\begin{bbook}[mr]
\bauthor{\bsnm{Bertoin},~\bfnm{Jean}\binits{J.}}
(\byear{1996}).
\btitle{L\'evy Processes}.
\bseries{Cambridge Tracts in Mathematics}
\bvolume{121}.
\bpublisher{Cambridge Univ. Press}, \baddress{Cambridge}.
\bid{mr={1406564}}
\bptok{imsref}%
\end{bbook}
\endbibitem

\bibitem{Blumenthal-Getoor-68}
\begin{bbook}[mr]
\bauthor{\bsnm{Blumenthal},~\bfnm{R.~M.}\binits{R.~M.}} \AND
  \bauthor{\bsnm{Getoor},~\bfnm{R.~K.}\binits{R.~K.}}
(\byear{1968}).
\btitle{Markov Processes and Potential Theory}.
\bseries{Pure and Applied Mathematics}
\bvolume{29}.
\bpublisher{Academic Press}, \baddress{New York}.
\bid{mr={0264757}}
\bptok{imsref}%
\end{bbook}
\endbibitem

\bibitem{Borodin-Salminen}
\begin{bbook}[mr]
\bauthor{\bsnm{Borodin},~\bfnm{Andrei~N.}\binits{A.~N.}} \AND
  \bauthor{\bsnm{Salminen},~\bfnm{Paavo}\binits{P.}}
(\byear{2002}).
\btitle{Handbook of {B}rownian Motion---Facts and Formulae}, \bedition{2nd} ed.
\bpublisher{Birkh\"auser}, \baddress{Basel}.
\bid{mr={1912205}}
\bptok{imsref}%
\end{bbook}
\endbibitem

\bibitem{Carmona-Petit-Yor-97}
\begin{bincollection}[mr]
\bauthor{\bsnm{Carmona},~\bfnm{Philippe}\binits{P.}},
  \bauthor{\bsnm{Petit},~\bfnm{Fr{\'e}d{\'e}rique}\binits{F.}} \AND
  \bauthor{\bsnm{Yor},~\bfnm{Marc}\binits{M.}}
(\byear{1997}).
\btitle{On the distribution and asymptotic results for exponential functionals
  of {L}\'evy processes}.
In \bbooktitle{Exponential Functionals and Principal Values Related to
  {B}rownian Motion}
\bpages{73--130}.
\bpublisher{Rev. Mat. Iberoam.}, \baddress{Madrid}.
\bid{mr={1648657}}
\bptok{imsref}%
\end{bincollection}
\endbibitem

\bibitem{DayKar}
\begin{barticle}[mr]
\bauthor{\bsnm{Dayanik},~\bfnm{Savas}\binits{S.}} \AND
  \bauthor{\bsnm{Karatzas},~\bfnm{Ioannis}\binits{I.}}
(\byear{2003}).
\btitle{On the optimal stopping problem for one-dimensional diffusions}.
\bjournal{Stochastic Process. Appl.}
\bvolume{107}
\bpages{173--212}.
\bid{doi={10.1016/S0304-4149(03)00076-0}, issn={0304-4149}, mr={1999788}}
\bptok{imsref}%
\end{barticle}
\endbibitem

\bibitem{Dynk}
\begin{barticle}[mr]
\bauthor{\bsnm{Dynkin},~\bfnm{E.~B.}\binits{E.~B.}}
(\byear{1963}).
\btitle{Optimal choice of the stopping moment of a {M}arkov process}.
\bjournal{Dokl. Akad. Nauk SSSR}
\bvolume{150}
\bpages{238--240}.
\bid{issn={0002-3264}, mr={0154329}}
\bptok{imsref}%
\end{barticle}
\endbibitem

\bibitem{Graversen-Peskir}
\begin{barticle}[mr]
\bauthor{\bsnm{Graversen},~\bfnm{S.~E.}\binits{S.~E.}} \AND
  \bauthor{\bsnm{Pe{\v{s}}kir},~\bfnm{G.}\binits{G.}}
(\byear{1997}).
\btitle{On {W}ald-type optimal stopping for {B}rownian motion}.
\bjournal{J. Appl. Probab.}
\bvolume{34}
\bpages{66--73}.
\bid{issn={0021-9002}, mr={1429055}}
\bptok{imsref}%
\end{barticle}
\endbibitem

\bibitem{hawkes}
\begin{barticle}[mr]
\bauthor{\bsnm{Hawkes},~\bfnm{John}\binits{J.}}
(\byear{1979}).
\btitle{Potential theory of {L}\'evy processes}.
\bjournal{Proc. Lond. Math. Soc. (3)}
\bvolume{38}
\bpages{335--352}.
\bid{doi={10.1112/plms/s3-38.2.335}, issn={0024-6115}, mr={0531166}}
\bptok{imsref}%
\end{barticle}
\endbibitem

\bibitem{MR2083932}
\begin{barticle}[mr]
\bauthor{\bsnm{Irle},~\bfnm{A.}\binits{A.}} \AND
  \bauthor{\bsnm{Paulsen},~\bfnm{V.}\binits{V.}}
(\byear{2004}).
\btitle{Solving problems of optimal stopping with linear costs of
  observations}.
\bjournal{Sequential Anal.}
\bvolume{23}
\bpages{297--316}.
\bid{doi={10.1081/SQA-200027048}, issn={0747-4946}, mr={2083932}}
\bptok{imsref}%
\end{barticle}
\endbibitem

\bibitem{ItoMckean}
\begin{bbook}[auto:STB|2011/08/30|13:45:12]
\bauthor{\bsnm{It{\^o}},~\bfnm{K.}\binits{K.}} \AND
  \bauthor{\bsnm{McKean},~\bfnm{H.~P.}\binits{H.~P.}}
(\byear{1996}).
\btitle{Diffusion Processes and Their Sample Paths}.
\bpublisher{Springer}, \baddress{Berlin}.
\bptok{imsref}%
\end{bbook}
\endbibitem

\bibitem{Khosh-Salm-Yor-06}
\begin{barticle}[mr]
\bauthor{\bsnm{Khoshnevisan},~\bfnm{Davar}\binits{D.}},
  \bauthor{\bsnm{Salminen},~\bfnm{Paavo}\binits{P.}} \AND
  \bauthor{\bsnm{Yor},~\bfnm{Marc}\binits{M.}}
(\byear{2006}).
\btitle{A note on a.s. finiteness of perpetual integral functionals of
  diffusions}.
\bjournal{Electron. Commun. Probab.}
\bvolume{11}
\bpages{108--117 (electronic)}.
\bid{issn={1083-589X}, mr={2231738}}
\bptok{imsref}%
\end{barticle}
\endbibitem

\bibitem{Kramkov-Morde-94}
\begin{barticle}[mr]
\bauthor{\bsnm{Kramkov},~\bfnm{D.~O.}\binits{D.~O.}} \AND
  \bauthor{\bsnm{Mordetski},~\bfnm{{\`E}.}\binits{{\`E}.}}
(\byear{1994}).
\btitle{An integral option}.
\bjournal{Teor. Veroyatn. Primen.}
\bvolume{39}
\bpages{201--211}.
\bid{doi={10.1137/1139007}, issn={0040-361X}, mr={1348195}}
\bptok{imsref}%
\end{barticle}
\endbibitem

\bibitem{Kyprianou-06}
\begin{bbook}[mr]
\bauthor{\bsnm{Kyprianou},~\bfnm{Andreas~E.}\binits{A.~E.}}
(\byear{2006}).
\btitle{Introductory Lectures on Fluctuations of {L}\'evy Processes with
  Applications}.
\bpublisher{Springer}, \baddress{Berlin}.
\bid{mr={2250061}}
\bptok{imsref}%
\end{bbook}
\endbibitem

\bibitem{kyprianoupistorius}
\begin{barticle}[mr]
\bauthor{\bsnm{Kyprianou},~\bfnm{A.~E.}\binits{A.~E.}} \AND
  \bauthor{\bsnm{Pistorius},~\bfnm{M.~R.}\binits{M.~R.}}
(\byear{2003}).
\btitle{Perpetual options and {C}anadization through fluctuation theory}.
\bjournal{Ann. Appl. Probab.}
\bvolume{13}
\bpages{1077--1098}.
\bid{doi={10.1214/aoap/1060202835}, issn={1050-5164}, mr={1994045}}
\bptok{imsref}%
\end{barticle}
\endbibitem

\bibitem{Lamperti-67}
\begin{bincollection}[mr]
\bauthor{\bsnm{Lamperti},~\bfnm{John}\binits{J.}}
(\byear{1967}).
\btitle{On random time substitutions and the {F}eller property}.
In \bbooktitle{Markov {P}rocesses and {P}otential {T}heory ({P}roc. {S}ympos.
  {M}ath. {R}es. {C}enter, {M}adison, {W}is., 1967)}
\bpages{87--101}.
\bpublisher{Wiley}, \baddress{New York}.
\bid{mr={0230370}}
\bptok{imsref}%
\end{bincollection}
\endbibitem

\bibitem{Lamperti-72}
\begin{barticle}[mr]
\bauthor{\bsnm{Lamperti},~\bfnm{John}\binits{J.}}
(\byear{1972}).
\btitle{Semi-stable {M}arkov processes. {I}}.
\bjournal{Z. Wahrsch. Verw. Gebiete}
\bvolume{22}
\bpages{205--225}.
\bid{mr={0307358}}
\bptok{imsref}%
\end{barticle}
\endbibitem

\bibitem{Lebedev-72}
\begin{bbook}[mr]
\bauthor{\bsnm{Lebedev},~\bfnm{N.~N.}\binits{N.~N.}}
(\byear{1972}).
\btitle{Special Functions and Their Applications}.
\bpublisher{Dover}, \baddress{New York}.
\bid{mr={0350075}}
\bptok{imsref}%
\end{bbook}
\endbibitem

\bibitem{Patie-ouq-06}
\begin{barticle}[mr]
\bauthor{\bsnm{Patie},~\bfnm{P.}\binits{P.}}
(\byear{2008}).
\btitle{{$q$}-invariant functions for some generalizations of the
  {O}rnstein--{U}hlenbeck semigroup}.
\bjournal{ALEA Lat. Am. J. Probab. Math. Stat.}
\bvolume{4}
\bpages{31--43}.
\bid{issn={1980-0436}, mr={2383732}}
\bptok{imsref}%
\end{barticle}
\endbibitem

\bibitem{Patie-06-poch}
\begin{barticle}[mr]
\bauthor{\bsnm{Patie},~\bfnm{Pierre}\binits{P.}}
(\byear{2009}).
\btitle{Exponential functional of a new family of {L}\'evy processes and
  self-similar continuous state branching processes with immigration}.
\bjournal{Bull. Sci. Math.}
\bvolume{133}
\bpages{355--382}.
\bid{doi={10.1016/j.bulsci.2008.10.001}, issn={0007-4497}, mr={2532690}}
\bptok{imsref}%
\end{barticle}
\endbibitem

\bibitem{Patie-06c}
\begin{barticle}[mr]
\bauthor{\bsnm{Patie},~\bfnm{Pierre}\binits{P.}}
(\byear{2009}).
\btitle{Infinite divisibility of solutions to some self-similar
  integro-differential equations and exponential functionals of {L}\'evy
  processes}.
\bjournal{Ann. Inst. Henri Poincar\'e Probab. Stat.}
\bvolume{45}
\bpages{667--684}.
\bptok{imsref}%
\end{barticle}
\endbibitem

\bibitem{patievigon}
\begin{bmisc}[auto:STB|2011/08/30|13:45:12]
\bauthor{\bsnm{Patie},~\bfnm{P.}\binits{P.}} \AND
  \bauthor{\bsnm{Vigon},~\bfnm{V.}\binits{V.}}
(\byear{2011}).
\bhowpublished{One-dimensional completely asymmetric
Markov processes. Unpublished manuscript, ULB.}
\bptok{imsref}%
\end{bmisc}
\endbibitem

\bibitem{PeskirShiryaev}
\begin{bbook}[mr]
\bauthor{\bsnm{Peskir},~\bfnm{Goran}\binits{G.}} \AND
  \bauthor{\bsnm{Shiryaev},~\bfnm{Albert}\binits{A.}}
(\byear{2006}).
\btitle{Optimal Stopping and Free-Boundary Problems}.
\bpublisher{Birkh\"auser}, \baddress{Basel}.
\bid{mr={2256030}}
\bptok{imsref}%
\end{bbook}
\endbibitem

\bibitem{Pitman-Yor-81}
\begin{bincollection}[mr]
\bauthor{\bsnm{Pitman},~\bfnm{Jim}\binits{J.}} \AND
  \bauthor{\bsnm{Yor},~\bfnm{Marc}\binits{M.}}
(\byear{1981}).
\btitle{Bessel processes and infinitely divisible laws}.
In \bbooktitle{Stochastic Integrals ({P}roc. {S}ympos., {U}niv. {D}urham,
  {D}urham, 1980)}.
\bseries{Lecture Notes in Math.}
\bvolume{851}
\bpages{285--370}.
\bpublisher{Springer}, \baddress{Berlin}.
\bid{mr={0620995}}
\bptok{imsref}%
\end{bincollection}
\endbibitem

\bibitem{Prabhakar-71}
\begin{barticle}[mr]
\bauthor{\bsnm{Prabhakar},~\bfnm{Tilak~Raj}\binits{T.~R.}}
(\byear{1971}).
\btitle{A singular integral equation with a generalized {M}ittag {L}effler
  function in the kernel}.
\bjournal{Yokohama Math. J.}
\bvolume{19}
\bpages{7--15}.
\bid{issn={0044-0523}, mr={0293349}}
\bptok{imsref}%
\end{barticle}
\endbibitem

\bibitem{Rivero-05}
\begin{barticle}[mr]
\bauthor{\bsnm{Rivero},~\bfnm{V{\'{\i}}ctor}\binits{V.}}
(\byear{2005}).
\btitle{Recurrent extensions of self-similar {M}arkov processes and
  {C}ram\'er's condition}.
\bjournal{Bernoulli}
\bvolume{11}
\bpages{471--509}.
\bid{doi={10.3150/bj/1120591185}, issn={1350-7265}, mr={2146891}}
\bptok{imsref}%
\end{barticle}
\endbibitem

\bibitem{Rogers-Williams-1}
\begin{bbook}[mr]
\bauthor{\bsnm{Rogers},~\bfnm{L.~C.~G.}\binits{L.~C.~G.}} \AND
  \bauthor{\bsnm{Williams},~\bfnm{David}\binits{D.}}
(\byear{2000}).
\btitle{Diffusions, {M}arkov Processes, and Martingales. {V}ol. 1}.
\bpublisher{Cambridge Univ. Press}, \baddress{Cambridge}.
\bnote{Reprint of the second (1994) edition}.
\bid{mr={1796539}}
\bptok{imsref}%
\end{bbook}
\endbibitem

\bibitem{Salminen}
\begin{barticle}[mr]
\bauthor{\bsnm{Salminen},~\bfnm{P.}\binits{P.}}
(\byear{1985}).
\btitle{Optimal stopping of one-dimensional diffusions}.
\bjournal{Math. Nachr.}
\bvolume{124}
\bpages{85--101}.
\bid{doi={10.1002/mana.19851240107}, issn={0025-584X}, mr={0827892}}
\bptok{imsref}%
\end{barticle}
\endbibitem

\bibitem{Shepp-67}
\begin{barticle}[mr]
\bauthor{\bsnm{Shepp},~\bfnm{L.~A.}\binits{L.~A.}}
(\byear{1967}).
\btitle{A first passage problem for the {W}iener process}.
\bjournal{Ann. Math. Statist.}
\bvolume{38}
\bpages{1912--1914}.
\bid{issn={0003-4851}, mr={0217879}}
\bptok{imsref}%
\end{barticle}
\endbibitem

\bibitem{shiryayevbook}
\begin{bbook}[mr]
\bauthor{\bsnm{Shiryayev},~\bfnm{A.~N.}\binits{A.~N.}}
(\byear{1978}).
\btitle{Optimal Stopping Rules}.
\bpublisher{Springer}, \baddress{New York}.
\bid{mr={0468067}}
\bptok{imsref}%
\end{bbook}
\endbibitem

\bibitem{Williams74}
\begin{barticle}[mr]
\bauthor{\bsnm{Williams},~\bfnm{David}\binits{D.}}
(\byear{1974}).
\btitle{Path decomposition and continuity of local time for one-dimensional
  diffusions. {I}}.
\bjournal{Proc. Lond. Math. Soc. (3)}
\bvolume{28}
\bpages{738--768}.
\bid{issn={0024-6115}, mr={0350881}}
\bptok{imsref}%
\end{barticle}
\endbibitem

\end{thebibliography}
\end{document}